\documentclass{amsart}
\usepackage{amsfonts}
\usepackage{amsmath,amsthm}
\usepackage{amsfonts,amssymb}
\usepackage{CJK}
\usepackage{enumerate}

\newtheorem{thm}{Theorem}[section]
\newtheorem{cor}[thm]{Corollary}
\newtheorem{lem}[thm]{Lemma}

\theoremstyle{remark}
\newtheorem{rem}[thm]{Remark}

\numberwithin{equation}{section}

\newcommand{\al}{\alpha}

\def\vz{\varepsilon}
\def\oz{\omega}
\def\lz{\lambda}

\def\dz{\delta}
\def\az{\alpha}
\def\gz{\gamma}

\def\({\Bigl(}
\def \){ \Bigr)}

\def\ga{\gamma}

\def\gaa{\boldsymbol{\ga}}

\def\x{{\bf x}}
\def\y{{\bf y}}

 \def\a{{\alpha}}

 \def\RR{{\mathbb R}}

\def\R{\mathbb{R}}
\def\h{{\bf h}}

\def\va{\varepsilon}

\def\y{{\bf y}}
\def\j{{\bf j}}

\def\va{\varepsilon}

\begin{document}
\def\RR{\mathbb{R}}
\def\Exp{\text{Exp}}
\def\FF{\mathcal{F}_\al}

\title[] {Tractability of non-homogeneous tensor product problems in the worst case setting}

\author[]{ Rong Guo} \address{ School of Mathematical Sciences, Capital Normal
University, Beijing 100048,
 China.}
\author[]{ Heping Wang} \address{ School of Mathematical Sciences, Capital Normal
University, Beijing 100048,
 China.}
\email{ wanghp@cnu.edu.cn.}

\keywords{ Tractability; Multivariate problems; Worst case
setting}

\subjclass[2010]{41A25, 41A63, 65D15, 65Y20}

\begin{abstract}
We study multivariate linear tensor product problems with some
special properties in the worst case setting. We  consider
algorithms that use finitely many continuous linear functionals.
We use a unified method to  investigate tractability of the above
multivariate problems, and obtain
 necessary and sufficient conditions  for strong
polynomial tractability, polynomial tractability, quasi-polynomial
tractability, uniformly weak tractability,  $(s,t)$-weak
tractability, and  weak tractability. Our results can apply to
multivariate approximation problems with kernels corresponding to
Euler kernels,  Wiener kernels, Korobov kernels, Gaussian kernels,
and analytic Korobov kernels.

\end{abstract}

\maketitle
\input amssym.def

\section{Introduction }
Recently, there has been an increasing interest in $d$-variate
computational problems $S=\{S_d\}_{d\in \Bbb N}$ with large or
even huge $d$. In order to solve these problems we usually
consider  algorithms using finitely many information operations.
For a given error threshold $\vz\in(0,1)$, the information
complexity $n(\vz, S_d)$ is defined to be the minimal number of
information operations for which the approximation error of some
algorithm is at most $\vz$. It is interesting  to study how the
information complexity $n(\varepsilon,S_d)$ depends on
$\varepsilon^{-1}$ and $d$. This is the subject of tractability.
Research on tractability of multivariate problems started in 1994
(see \cite{W}) and they quickly gained the popularity.
 Nowadays,  tractability form an area of an intensive
research with wide scope of unexpected results and open problems.
The fundamental references about tractability are books \cite{NW1,
NW2, NW3}. There are new papers on the subject constantly coming
out (see e.g. \cite{CW1, CWZ,  DKPW, IKPW, LX1, S1, SiW, SW}).

 Various notions of
tractability have been studied recently for many multivariate
problems.  We  recall some of the basic tractability notions (see
\cite{NW1, NW2, NW3, S1, SiW}).
 Let
$S= \{S_d\}_{d\in\Bbb N}$. We say the problem $S$ is

$\bullet$ {\it strongly polynomially tractable (SPT)} iff there exist
non-negative numbers $C$ and $p$ such that
\begin{equation}\label{1.1} n(\varepsilon,S_d)\leq C(\varepsilon^{-1})^p \ \ {\rm for \ \ all} \ \ d\in\Bbb N,\ \ \varepsilon\in (0,1).\end{equation}
The exponent $p^*$ of SPT is defined as the infimum of $p$ for
which \eqref{1.1} holds;

$\bullet$ {\it polynomially tractable (PT)} iff there exist non-negative
numbers $C$, $p$ and $q$ such that
$$ n(\varepsilon,S_d)\leq Cd^q(\varepsilon^{-1})^p \ \ {\rm for \ \ all  }\ \ d\in\Bbb N,\ \ \varepsilon\in (0,1);$$

$\bullet$ {\it quasi-polynomially tractable (QPT)} iff there exist
positive numbers $C$ and $t$ such that
\begin{equation}\label{1.2}  n(\varepsilon,S_d)\leq C\exp(t(1+\ln d)(1+\ln\varepsilon^{-1})) \ \ {\rm for \ \ all  }\ \ d\in  \Bbb N,\ \ \varepsilon\in (0,1).\end{equation}
The exponent $t^*$ of QPT  is defined as the infimum of $t$ for
which \eqref{1.2} holds;

$\bullet$ {\it uniformly weakly tractable (UWT)} iff for all $\alpha,
\beta >0$
$$\lim _{\varepsilon^{-1}+d \rightarrow \infty}\frac{\ln n(\varepsilon,S_d)}{(\varepsilon^{-1})^\alpha+d^\beta}=0;$$

$\bullet$ {\it weakly tractable (WT)} iff
$$\lim _{\varepsilon^{-1}+d \rightarrow \infty}\frac{\ln n(\varepsilon,S_d)}{\varepsilon^{-1}+d}=0;$$

$\bullet$ {\it $(s,t)$-weakly tractable ($(s,t)$-WT)} for fixed
positive $s$ and $t$ iff
\begin{equation*}
\lim_{\varepsilon ^{-1}+d\rightarrow \infty }\frac{\ln n(\va
,S_d)}{(\va ^{-1})^{s }+d^{t }}=0.
\end{equation*}

Clearly, $(1,1)$-WT is the same as WT. If the problem $S$ is not
WT, then $S$ is called  intractable. We say that the problem $S$
suffers from {the curse of dimensionality} if there exist positive
numbers $C, \,\va _0,\, \a $ such that for all $0<\va\leq \va
_{0}$ and infinitely many $d\in \Bbb N$,
\begin{equation*}
n(\va ,S_d)\geq C(1+\a )^{d}.
\end{equation*}

This paper is devoted to studying  tractability of non-homogeneous
tensor product  problems   with  some special properties  in the
worst case setting. Such approximation problems were investigated
in \cite{PW, S2} for Korobov kernels, in \cite{FHW1} for Gaussian
kernels, and in \cite{IKPW, KPW, LX1} for analytic Korobov
kernels. In the average case setting, there are many papers
devoted to discussing tractability of non-homogeneous tensor
product  problems with covariance kernels being  a product of
univariate kernels and satisfying some special properties (see
\cite{CW1, CWZ, FHW2, LPW1, LPW2, LX, S2, S3, X1, X2}).

In this paper, we use a unified method to study tractability of
the above multivariate problems in the worst case setting, and
obtain
 necessary and sufficient conditions  for strong
polynomial tractability, polynomial tractability, quasi-polynomial
tractability, uniformly weak tractability,  $(s,t)$-weak
tractability, and  weak tractability. Our results can apply to
multivariate approximation problems with kernels corresponding to
Euler kernels,  Wiener kernels, Korobov kernels, Gaussian kernels,
and analytic Korobov kernels. For the first three kernels the
obtained results are new.

The paper is organized as follows. In Section 2 we give
preliminaries about
 non-homogeneous tensor product problems  in the
worst case setting and  present the main results, i.e.,  Theorem
2.1.
   Section 3 is
devoted to proving Theorem 2.1.
  In  Section 4, we give the applications of Theorem 2.1 to  the  problems
  with
kernels corresponding to Euler kernels,  Wiener kernels, Korobov
kernels, Gaussian kernels, and analytic Korobov kernels.

\section{Preliminaries and main results}

In this section, we recall the definition of non-homogeneous
linear tensor product problem in the worst case setting (see
\cite{PW, S2}).

A linear tensor product problem in the worst case setting is a sequence
of linear operators
$$S=\{S_d\}_{d\in \Bbb N}$$
such that for every $j\in \Bbb N$ there exists a separable Hilbert
space $\mathcal{H}_j$, a Hilbert space $\mathcal{G}_j$ and a
continuous linear operator
$\mathcal{S}_j:\mathcal{H}_j\to\mathcal{G}_j$ such that
$$S_d=\bigotimes_{j=1}^d \mathcal{S}_j: H_d\to G_d,$$
where $H_d=\bigotimes_{j=1}^d\mathcal{H}_j$ and
 $G_d=\bigotimes_{j=1}^d\mathcal{G}_j$ for every $d\in \Bbb
 N$.

 If $\mathcal{H}_j=\mathcal{H}_1$, $\mathcal{G}_j=\mathcal{G}_1$ and
 $\mathcal{S}_j=\mathcal{S}_1$ for all $j\in \Bbb N$,
 then the linear tensor product problem $S$ is called homogeneous.

 Without loss of generality we assume that all operators
$\mathcal{S}_k,\ k\in \Bbb N$ are compact. Then the operators
$$\mathcal{W}_k=\mathcal{S}_k^*\mathcal{S}_k:\mathcal{H}_k\to\mathcal{H}_k\
\  k\in \Bbb N$$ are compact,  self-adjoint, and non-negative
definite, where $\mathcal{S}_k^*$ are the adjoint operators of
$\mathcal{S}_k$.  Let $\{(\lambda(k,j),\eta(k,j))\}_{j\in  \Bbb
N}$ denote the  eigenpairs of $\mathcal{W}_k$ satisfying
$$
\mathcal W_k(\eta(k,j))=\lambda(k,j)\eta(k,j) \ \ {\rm with } \ \
\lambda(k,1)\geq \lambda(k,2)\geq \cdots\geq 0, \
\eta(k,j)\in\mathcal H_k .$$ Then the eigenpairs of the operators
$$W_d=S_d^*S_d=\bigotimes_{j=1}^d\mathcal{W}_j,\ d\in\Bbb N$$ are given by
$$\big\{(\lz _{d,\j},\eta _{d,\j})\big\}_{\j=(j_1,j_2,\dots,j_d)\in \Bbb N^d},$$
where
$$\lz_{d,\j}=\prod _{k=1}^d\lz (k,j_k),\quad \quad \eta_{d,\j}=\prod_{k=1}^d\eta(k,j_k).$$

Let the sequence $\{\lz _{d,j}\}_{j\in \Bbb N}$ be the
nonincreasing rearrangement of $\{\lz_{d,\j}\}_{\j\in \Bbb N^d}$.
Then we have
\begin{equation}\label{2.1}
\sum_{j=1}^\infty\lz^{\tau}_{d,j}=\prod_{k=1}^d\sum_{j=1}^\infty
\lz(k,j)^\tau,\quad \text{for any}\quad \tau>0.
\end{equation}
Note that the above sums are not always finite for any $\tau>0$.

  We approximate $S_df$ by algorithms
$A_{n,d}f$ of the form
\begin{equation}\label{2.2}A_{n,d}f=\phi
_{n,d}(L_1(f),L_2(f),\dots,L_n(f)),\end{equation} where
$L_1,L_2,\dots,L_n$ are continuous linear functionals on $H_d$,
 and $\phi _{n,d}:\;\Bbb R^n\to
G_d$ is an arbitrary measurable mapping. The worst case
approximation error for the algorithm $A_{n,d}$  is defined as
 $$e(A_{n,d})=\sup_{\|f\|_{H_d}\le1}
\|S_d\,f-A_{n,d}f\|_{G_d}.$$ The $n$th minimal worst case error,
for $n\ge 1$, is defined by
$$e(n,S_d)=\inf_{A_{n,d}}e(A_{n,d}),$$
where the infimum is taken over all algorithms of the form
\eqref{2.2}. According to \cite{NW1},  the  $n$th minimal worst
case error is given by
$$e(n,S_d)=\lz_{d,n+1}^{1/2},$$and
is achieved by the $n$th optimal algorithm
\begin{equation*}
A_{n,d}^{*}f=\sum_{j=1}^n  \big\langle f,\eta _{d,j} \big \rangle
_{H_d}S_d\,\eta _{d,j}.\end{equation*}

For $n=0$, we use $A_{0,d}=0$. We remark that  the so-called
initial error $e(0,S_d )$ is  defined by
$$e(0,S_d )=\sup_{\|f\|_{H_d}\le1}\|S_d\,
f\|_{G_d}=\|S_d\|=\lz_{d,1}^{1/2}.$$

The information complexity for $S_d$ can be studied using either
the absolute error criterion (ABS), or the normalized error
criterion (NOR). We define the information complexity $n^X(\va
,S_d)$ for $X\in \{ {\rm ABS,\, NOR}\}$ as
\begin{equation*}
n^{ X}(\va ,S_d)=\min\{n:\,e(n,S_d)\leq \va CRI_d\},
\end{equation*}where
\begin{equation*}
CRI_d=\left\{\begin{matrix}
 & 1, \; \quad\qquad\text{ for X=ABS,} \\
 &e(0,S_d),\quad \text{ for X=NOR.}
\end{matrix}\right.
\end{equation*}

In this paper we consider a special class of non-homogeneous
linear tensor product problems
 $S=\{S_d\}_{d\in \Bbb N}$ in the worst case setting. In the average case setting the similar non-homogeneous
linear tensor product problems were investigated in \cite{CWZ}.
Assume that the eigenvalues
$$\big\{\prod_{k=1}^d\lambda(k,j_k)\big\}_{(j_1,j_2,\cdots,j_d)\in  \Bbb N^d}$$
of the operator $W_d=S_d^*S_d$ of the problem  $S=\{S_d\}_{d\in
\Bbb N}$ satisfy the following three conditions:

(1) $\lambda(k,1)=1, k\in  \Bbb N$;

(2) the sequence $\{h_k\}_{k\in \Bbb N}$ is  nonincreasing, where
$h_k=\frac{\lambda(k,2)}{\lambda(k,1)}\in(0,1]$;

(3) there exist a positive constant $\tau>0$ such that
\begin{equation}\label{2.3}\sup_{{k\in  \Bbb N}}H(k,\tau)=
H(1,\tau)\ \ {\rm and}\ \
M_\tau=H(1,\tau)<\infty,\end{equation}where
$$H(k,\tau)=\sum\limits_{j=2}^\infty\Big(\frac{\lambda(k,j)}{\lambda(k,2)}\Big)^\tau.$$
Let $$\tau_0=\inf\big\{\tau\ \big|\ \tau \ \  {\rm satisfies}\ \
\eqref{2.3}\big\}.$$ Then we say that the problem $S=\{S_d\}_{d\in
\Bbb N}$ has Property (P).

We are ready to present the main result of this paper.

\begin{thm}\label{thm2.1} Let $S=\big\{S_d\big\}_{d\in\Bbb N}$ be a non-homogeneous tensor
product problem with Property (P) in the worst case setting.  Then
for the absolute error criterion or the normalized error
criterion, we have

(i) $S$ is strongly polynomially tractable iff
\begin{equation} A_*=\liminf\limits_ {d\to\infty}\frac{\ln h_d^{-1}}{\ln d}>0,\label{2.4}\end{equation}
 and the exponent of SPT is \begin{equation}\label{2.5} p^*=\max\big\{\frac{2}{A_*},2\tau_0\big\}.\end{equation}

(ii) $S$ is strongly polynomially tractable iff it is polynomially
tractable.

\vskip 2mm

(iii)  $S$ is quasi-polynomially tractable iff
\begin{equation}\label{2.6}B:=\lim\limits_{d\to\infty}\ln h_d^{-1}>0.\end{equation}This is equivalent to that $h_k\not\equiv1$. Furthermore, the exponent of QPT is
\begin{equation}\label{2.7}  t^*=\max\big\{\frac{2}{B},2\tau_0\big\}.\end{equation}

(iv) $S$ is quasi-polynomially tractable iff $S$ is uniformly
weakly tractable, iff $S$ is $(s,t)$-weakly tractable with $s>0$
and $t\in (0,1]$, and  iff $S$ is weakly tractable.

\vskip 2mm

(v) $S$ is $(s,t)$-weakly tractable with $s>0$ and $t>1$.

\vskip 2mm

 (vi)  $S$ suffers from the curse of dimensionality iff
$h_k\equiv1$.

\end{thm}

\begin{rem} Let $S=\big\{S_d\big\}_{d\in\Bbb N}$ be a non-homogeneous tensor
product problem. If the eigenvalues of the  operator
$W_d=S_d^*S_d$  satisfy Conditions (2) and (3) of Property (P),
then for NOR, Theorem \ref{thm2.1} holds.

Indeed, let $\tilde S=\big\{\tilde S_d\big\}_{d\in\Bbb N}$ be the
non-homogeneous tensor product problem which  the eigenvalues
$\big\{\prod_{k=1}^d\tilde
\lz(k,j_k)\big\}_{(j_1,j_2,\dots,j_d)\in \Bbb N^d}$ of the
corresponding  operator $\tilde W_d=\tilde S_d^* \tilde S_d$ of
  satisfy $$ \tilde
\lz(k,j)=\frac {\lz(k,j)}{ \lz(k,1)},\ \  j\in\Bbb N,\
k=1,\dots,d.$$Then $\tilde S$ has Property (P) with the same $h_k$
and Theorem \ref{thm2.1} is applicable. Also for NOR, the problems
$S$ and $\tilde S$ have the same tractability. Hence, for NOR,
Theorem \ref{thm2.1} holds.

\end{rem}

In order to prove Theorem \ref{thm2.1} we need the following
lemmas.

\begin{lem} (See \cite[Theorem 5.2]{NW1}.)  Consider the non-zero problem $S=\{S_d\}_{d\in \Bbb N}$ in the worst case
setting, where $S_d: H_d \to G_d$ is a compact linear operator
between two Hilbert spaces $H_d$ and $G_d$. Assume that  the
eigenvalues $\{\lz_{d,j}\}_{j\in \Bbb N}$ of the operator
$W_d=S_d^*S_d$ satisfying $$\lz_{d,1}\ge \lz_{d,2}\ge \cdots\ge 0
. $$ Then  for the normalized error criterion, we have

(a) $S$ is polynomially tractable iff there exist two constants
$q\ge 0$ and $\tau
> 0$ such that  \begin{equation}\label{2.20} C_{\tau,q}:=\sup_{d\in  \Bbb N}\Big(\sum\limits_{j=1}^\infty
\Big(\frac{\lambda_{d,j}}{\lambda_{d,1}}\Big)^\tau\Big)^\frac{1}{\tau}d^{-q}<\infty.\end{equation}

(b) $S$ is strongly polynomially tractable iff \eqref{2.20} holds
with $q=0$. If so then the exponent of SPT is
$$p^*=\inf\big\{2\tau\ \big|\ \tau \ {\rm satisfies}\
\eqref{2.20}\ {\rm with}\ q=0\big\}.$$\end{lem}

\begin{lem} (See \cite[Theorem 23.2]{NW3}.) Consider the  problem $S=\{S_d\}_{d\in \Bbb N}$ in the worst case
setting as in Lemma 2.3. Then for the normalized error criterion,
$S$ is quasi-polynomially tractable iff there exists a $\tau
> 0$ such that
\begin{equation}\label{2.21}B_\tau:=\sup_{d\in  \Bbb N}
\Big(\sum\limits_{j=1}^\infty
\Big(\frac{\lambda_{d,j}}{\lambda_{d,1}}\Big)^{\tau(1+\ln
d)}\Big)^\frac{1}{\tau}d^{-2}<\infty.\end{equation}If so then the
exponent of QPT is
$$t^*=\inf\big\{2\tau\ \big|\ \tau \ {\rm satisfies}\
\eqref{2.21}\big\}.$$\end{lem}

\begin{lem} Consider the  problem $S=\{S_d\}_{d\in \Bbb N}$ in the worst case
setting as in Lemma 2.3. For any fixed $t>0$ if there exists a
positive $\tau$ such that
$$\lim\limits_{d\to\infty}d^{-t}\ln\Big( \sum\limits_{j=1}^\infty\Big(\frac{\lambda_{d,j}}{\lambda_{d,1}}\Big)^\tau\Big)=0,$$
then   $S$ is $(s,t)$-weakly tractable for any $s>0$ and this $t$
for the normalized error criterion.
\end{lem}
Lemma 2.5 is likely not new, however we cannot find its proof. For
readers' convenience we give its proof.

\begin{proof} Without loss of generality we assume that
$\lz_{d,1}=1$. Since
$$ n\ \lambda_{d,n}^\tau\leq\sum\limits_{k=1}^n\lambda_{d,k}^\tau\leq\sum\limits_{k=1}^\infty\lambda_{d,k}^\tau,$$
we get
$$\lambda_{d,n}\leq n^{-\frac{1}{\tau}}\Big(\sum\limits_{k=1}^\infty\lambda_{d,k}^\tau\Big)^{\frac 1\tau}.$$
We note that $\lambda_{d,n}\leq\varepsilon^2$ holds for
$$n=\Big\lceil\sum\limits_{k=1}^\infty\lambda_{d,k}^\tau \ \varepsilon^{-2\tau}\Big\rceil.$$
It follows that
$$n(\varepsilon,S_d)=\min\{n|\lambda_{d,n+1}\leq\varepsilon^2\}\leq\Big\lceil\sum\limits_{k=1}^\infty\lambda_{d,k}^\tau \ \varepsilon^{-2\tau}\Big\rceil.$$
Since $\lceil x\rceil\leq x+1 \leq2\max\{1,x\}$ and $\ln\max\{x,1\}\leq\max\{0,\ln x\}=(\ln x)_+$, we have
$$\frac{\ln n(\varepsilon,S_d)}{\varepsilon^{-s}+d^t}
\leq\frac{\ln2+\Big(2\tau\ln\varepsilon^{-1}+ln\Big(\sum\limits_{k=1}^\infty\lambda_{d,k}^\tau\Big)\Big)_+}{\varepsilon^{-s}+d^t},$$
which goes to $0$ as $\varepsilon^{-1}+d$ tends to infinity. Hence $S$ is $(s,t)$-weakly tractable with $s>0$ and this $t$.

\end{proof}

\section{Proof of Theorem \ref{thm2.1}}

\noindent{\it \textbf{Proof of Theorem \ref{thm2.1}}.} \vskip 2mm

 Since $\lambda_{d,1}=\prod_{k=1}^d\lz(k,1)=1 $ we know that tractability for $S=\big\{S_d\big\}_{d\in\Bbb N}$ is
the same for both the absolute and normalized error criteria.

(i).  If $S$ is strongly polynomially tractable with the exponent
$p^*$ of SPT, then by Lemma 2.3 (b) we know that for any fixed
$x>{p^*}/{2}$, we have
\begin{equation*} C_x:=\sup_{d\in  \Bbb N}\Big(\sum\limits_{j=1}^\infty \lambda_{d,j}^x\Big)^\frac{1}{x}<\infty.\end{equation*}
It follows from \eqref{2.1} and  \eqref{2.3} that for any $ d\in
\Bbb N$,
\begin{align}\sum\limits_{j=1}^\infty \lambda_{d,j}^x & =\prod\limits_{k=1}^d\sum\limits_{j=1}^\infty \lambda(k,j)^x
                                     =\prod\limits_{k=1}^d(1+\sum\limits_{j=2}^\infty\lambda(k,j)^x)
                                     \notag \\ &=\prod\limits_{k=1}^d(1+H(k,x)h_k^x)\le C_x^x<\infty.\label{3.1}\end{align}
 This means that for $1\le k\le d$, $H(k,x)<\infty$. This means that $$x\ge
 \tau_0.$$
Clearly for any $x>{p^*}/{2}$, we have
\begin{equation*}H(k,x)=1+\sum\limits_{j=3}^\infty\Big(\frac{\lambda(k,j)}{\lambda(k,2)}\Big)^x\ge 1.\end{equation*}
 It follows from \eqref{3.1} that
\begin{equation}\label{3.2}\prod\limits_{k=1}^d(1+h_k^x)\leq \prod\limits_{k=1}^d(1+H(k,x)h_k^x)\le C_x^x. \end{equation}
Using the inequality $\ln(1+x)\geq x\ln2,\ x\in[0,1]$,  Condition
(2), and \eqref{3.2},  we  get  that
\begin{equation}\label{3.2-1}d\ln2h_d^x\leq \ln2\sum\limits_{k=1}^d h_k^x\leq \sum\limits_{k=1}^d\ln(1+h_k^x)=\ln\Big(\prod\limits_{k=1}^d(1+h_k^x)
\Big)\leq \ln C_x^x,\end{equation} which means that $$\ln
d+\ln(\ln2)-x\ln h_d^{-1}\le \ln(\ln C_x^x). $$It follows that
\begin{equation}\label{3.3}\frac{\ln h_d^{-1}}{\ln d}\geq\frac{\ln d+\ln(\ln2)-\ln(\ln C_x^x)}{x\ln d}.\end{equation}
Letting $d\rightarrow\infty $ in \eqref{3.3}, we get that
\begin{equation}\label{3.4}A_*=\liminf\limits_ {d\rightarrow\infty}\frac{\ln h_d^{-1}}{\ln d}\geq\frac{1}{x}>0,\end{equation}
i.e., \eqref{3.1} holds. Since \eqref{3.4} and the inequality
$x\ge \tau_0$ hold for any $x> p^*/2$ we obtain that
\begin{equation}\label{3.5}p^*\geq \max\big\{\frac{2}{A_*},2\tau_0\big\}.\end{equation}

On the other hand, assume that  \eqref{2.4} holds. From
\eqref{2.3} we  know that for $x>\tau_0$,
 \begin{equation}\label{3.6}1\leq \sup_{{k\in  \Bbb N}}H(k,x)=\sup_{{k\in  \Bbb N}}\sum\limits_{j=2}^\infty\Big(\frac{\lambda(k,j)}{\lambda(k,2)}\Big)^x\leq M_x<\infty.\end{equation}
It follows from \eqref{3.1},  \eqref{3.6}, and the inequality
$\ln(1+x)\le x\ (x\ge0)$ that for any
$x>\max\big\{\frac{1}{A_*},\tau_0\big\}$ and $d\in\Bbb N$,
\begin{equation}\label{3.7}\sum\limits_{j=1}^\infty \lambda_{d,j}^x\leq \prod\limits_{k=1}^d(1+M_x
h_k^x)=\exp\Big(\sum\limits_{k=1}^d \ln (1+M_xh_k^x)\Big) \leq
\exp\Big(M_x \sum\limits_{k=1}^\infty h_k^x\Big) .\end{equation}
For any $x>\max\big\{\frac{1}{A_*},\tau_0\big\}$, by \eqref{2.4}
we have
$$ h_k^x=k^{-\frac{x \ln h_k^{-1}}{\ln k}}\ \ {\rm and}\ \
\liminf\limits_{k\to\infty}\frac{x\ln h_k^{-1}}{\ln
k}=xA_*>1,$$which implies that
 $$D_x:=\sum\limits_{k=1}^\infty h_k^x<\infty.$$ From \eqref{3.7} we  get that
 \begin{equation*}\sup_{d\in\Bbb N}\sum\limits_{j=1}^\infty\lambda_{d,j}^x\leq \exp(D_xM_x)<\infty.\end{equation*}
 By Lemma 2.3 (b) we obtain that  $S$ is strongly polynomially
 tractable, and
  the exponent of SPT satisfies
 \begin{equation}\label{3.8}p^*\leq \max\big\{\frac{2}{A_*},2\tau_0\big\}.\end{equation}
 From \eqref{3.5} and \eqref{3.8} we obtain \eqref{2.5}. The proof
 of (i) is complete.

 \vskip 3mm

 (ii). We know that SPT can deduce PT.  So it suffices to show that  PT implies SPT.  Assume that $S$ is PT.  By Lemma 2.3 (a)  there exist  numbers
  $x> 0,\ q\ge 0$  such that
  \begin{equation*} C_{x,q}:=\sup_{d\in  \Bbb N}\Big(\sum\limits_{j=1}^\infty \lambda_{d,j}^x\Big)^\frac{1}{x}d^{-q}<\infty.\end{equation*}
It follows from \eqref{2.1} that
\begin{equation}\label{3.10} \sum\limits_{k=1}^d\ln(1+h_k^x)\le \sum\limits_{k=1}^d\ln (1+H(k,x)h_k^x)=\sum\limits_{j=1}^\infty \lambda_{d,j}^x\leq
\ln(C_{x,q}^x d^{qx}).\end{equation}
  Similar to the proof of \eqref{3.2-1}, by \eqref{3.10} we  get
\begin{equation*}d\ln2h_d^x\leq \ln2\sum\limits_{k=1}^d h_k^x\leq \sum\limits_{k=1}^d\ln(1+h_k^x)\leq \ln (C_{x,q}^x)+qx\ln d,\end{equation*}
which implies that
$$\ln
d+\ln(\ln2)-x\ln h_d^{-1}\le \ln(\ln C_{x,q}^x+qx\ln d).
$$It follows that
\begin{equation}\label{3.11}\frac{\ln h_d^{-1}}{\ln d}\geq\frac{\ln d+\ln(\ln2)-\ln(\ln C_{x,q}^x+qx\ln
d)}{x\ln d}.\end{equation} Letting $d\rightarrow\infty $ in
\eqref{3.11}, we get that
  \begin{equation*} \liminf\limits_{d\to\infty}\frac{\ln{h_d^{-1}}}{\ln d}\geq \frac{1}{x}> 0.\end{equation*}
  This implies that \eqref{2.4} holds and hence by the proved (i) we obtain that $S$ is strongly polynomially tractable. This completes the proof of (ii).

   \vskip 3mm

(iii).  Assume that $S$ is quasi-polynomially tractable with the
exponent $t^*$ of QPT. Then by Lemma 2.4 we know that for any
fixed $ x>\frac{t^*}{2}$, we have
 \begin{equation*}B_x:=\sup_{d\in  \Bbb N}\Big(\sum\limits_{j=1}^\infty \lambda_{d,j}^{x(1+\ln d)}\Big)^\frac{1}{x}d^{-2}<\infty.\end{equation*}
It follows from \eqref{2.1} that
\begin{align} \sum\limits_{k=1}^d\ln\big(1+h_k^{x(1+\ln d)}\big)&\le \sum\limits_{k=1}^d\ln \Big(1+H(k,{x(1+\ln d)})\,h_k^{x(1+\ln d)}\Big)
\notag \\ &=\ln\Big(\sum\limits_{j=1}^\infty
\lambda_{d,j}^{x(1+\ln d)}\Big)\leq \ln(B_x^x)+ 2x\ln
d.\label{3.12} \end{align}If $d=1$, then the above inequality
gives that $$\sum\limits_{j=1}^\infty \lambda_{1,j}^{x(1+\ln
1)}=\sum\limits_{j=1}^\infty \lz(1,j)^x=1+H(1,x)h_1^x<\infty.
$$ This means that $x\ge \tau_0$. Using \eqref{3.12}, the inequality $\ln(1+x)\geq x\ln2 \
(x\in[0,1])$, and the monotonicity of the sequence of $\{h_k\}$,
we get
 \begin{equation*}\ln2\, d h_d^{x(1+\ln d)}\le \ln2  \sum\limits_{k=1}^d h_k^{x(1+\ln d)}\leq \sum\limits_{k=1}^d\ln\big(1+h_k^{x(1+\ln d)}\big)
 \leq \ln(B_x^x)+ 2x\ln
d, \end{equation*}which implies that
$$\ln \ln2+\ln d- x(1+\ln d)\ln h_d^{-1}\le \ln\big(\ln(B_x^x)+ 2x\ln
d\big).$$ This yields that
\begin{equation}\label{3.13} \ln h_d^{-1} \geq \frac{\ln d+\ln\ln2-\ln\big(\ln(B_x^x)+ 2x\ln
d\big)}{x(1+\ln d)}.\end{equation} Letting $d\rightarrow\infty$ in
\eqref{3.13} and noting the monotonicity of $\{h_k\}$, we get
\begin{equation*}B:=\lim\limits_ {d\rightarrow\infty}\ln h_d^{-1}\geq \frac{1}{x}>0 .\end{equation*}Since the
inequalities $B\ge 1/x$ and  $x\ge \tau_0$ hold for any $x> t^*/2$
we obtain that
\begin{equation}\label{3.14}t^*\geq \max\big\{\frac{2}{B},2\tau_0\big\}.\end{equation}

On the other hand, assume that  \eqref{2.6} holds. From
\eqref{2.3} we  know that for $x>\tau_0$,
 \begin{equation}\label{3.15}1\leq \sup_{{k\in  \Bbb N}}H(k,{x(1+\ln d)})\le \sup_{{k\in  \Bbb N}}H(k,{x})\leq M_x<\infty.\end{equation}
It follows from \eqref{3.12},  \eqref{3.15}, and the inequality
$\ln(1+x)\le x\ (x\ge0)$ that for any
$x>\max\big\{\frac{1}{B},\tau_0\big\}$ and $d\in\Bbb N$,
\begin{align}\sum\limits_{j=1}^\infty\lambda_{d,j}^{x(1+\ln d)}&\leq \prod\limits_{k=1}^d(1+M_x
h_k^{x(1+\ln d)})  =\exp\Big(\sum\limits_{k=1}^d \ln
(1+M_xh_k^{x(1+\ln d)})\Big)\notag \\&\leq \exp\Big(M_x
\sum\limits_{k=1}^d h_k^{x(1+\ln d)}\Big) \leq \exp\Big(M_x
\sum\limits_{k=1}^d h_k^{x(1+\ln k)}\Big)\notag\\ & \leq
\exp\Big(M_x \sum\limits_{k=1}^\infty h_k^{x(1+\ln
k)}\Big).\label{3.16}\end{align} For any
$x>\max\big\{\frac{1}{B},\tau_0\big\}$, by \eqref{2.6} we have
$$ h_k^{x(1+\ln k)}= k^{-\frac{{x(1+\ln k)}\ln h_k^{-1}}{\ln k}}\ \ {\rm and}\ \
\liminf\limits_ {k\rightarrow\infty}\frac{{x(1+\ln k)}\ln h_k^{-1}}{\ln k} = xB >1,$$
which implies that
 $$E_x:=\sum\limits_{k=1}^\infty h_k^{x(1+\ln k)}<\infty.$$ From \eqref{3.16} we  get that
 \begin{equation*}\sup_{d\in\Bbb N}\sum\limits_{j=1}^\infty \lambda_{d,j}^{x(1+\ln d)}\leq \exp(E_xM_x)<\infty.\end{equation*}
 By Lemma 2.4 we obtain that  $S$ is quasi-polynomially tractable, and
  the exponent of QPT satisfies
 \begin{equation}\label{3.17}t^*\leq \max\big\{\frac{2}{B},2\tau_0\big\}.\end{equation}
 From \eqref{3.14} and \eqref{3.17} we obtain \eqref{2.7}. This completes the  proof
 of (iii).

\vskip 3mm

(iv).  Since QPT $\Rightarrow$ UWT $\Rightarrow$ $(s,t)$-WT (or
WT), it suffices to show that $(s,t)$-WT $\Rightarrow$ QPT with
$t\in(0,1]$ and $s>0$. Assume that $S$ is $(s,t)$-weakly tractable
with $t\in(0,1]$ and $s>0$. First we show that $h_k\not\equiv 1$.
 If $h_k\equiv 1$, then $e(n,S_d)=1$ for $1\le n \le
2^d$. This means that
$$n(\vz,S_d)\ge 2^d-1, \ \ \vz \in (0,1),$$and
$$\lim_{d\to\infty}\frac{\ln n(1/2,S_d)}{2^s+d^t} \ge
\lim_{d\to\infty}d^{1-t}\ln2>0,\ \ t\in(0,1].$$ Hence $S$ suffers
from the curse of dimensionality and  is not $(s,t)$-weakly
tractable with $t\in(0,1]$ and $s>0$. This leads to a
contradiction. Hence $h_k\not\equiv 1$. By the monotonicity of
$\{h_k\}$ we get
$$B=\lim_{k\to\infty}\ln h_k^{-1}>0.$$ By the proved (iii) we obtain that $S$
is quasi-polynomially tractable. (iv) is proved.

\vskip 3mm

(v).
 It follows from  \eqref{3.7} that for $x>\tau_0$,
 $$\ln\Big(\sum\limits_{j=1}^\infty
 \lambda_{d,j}^x\Big)\le  {M_x}\sum\limits_{k=1}^dh_k^x\le  {M_x}d. $$
 We have for $t>1$,
 $$\lim_{d\to\infty}d^{-t}\ln\Big(\sum\limits_{j=1}^\infty
 \lambda_{d,j}^x\Big)\le \lim_{d\to\infty}
 {M_x}d^{1-t}=0.$$
 By lemma 2.5 we get that $S$ is $(s,t)$-weakly tractable with $s>0$ and $t>1$.
 This complete the proof of (v).

\vskip 3mm

 (vi). We have proved in (iv) that if $h_k\equiv 1$ then $S$ suffers
from the curse of dimensionality, and  if $h_k\not\equiv 1$ then
 $S$
is quasi-polynomially tractable and does not suffer from the curse
of dimensionality. This completes the proof of (vi).

\vskip 1mm

 The proof of Theorem \ref{thm2.1} is finished.

 \begin{rem} Let $S=\big\{S_d\big\}_{d\in\Bbb N}$ be a non-homogeneous tensor
product problem. If the eigenvalues of the  operator
$W_d=S_d^*S_d$  satisfy Condition (3) of Property (P) and the
following  Condition $(2)'$:
 there exist a nonincreasing positive
sequence $\big\{f_k\big\}_{k\in\Bbb N}$ and two positive constants
$A_1,\ A_2$ such that
 such that for all $k\in\Bbb N$, we have
$$A_1f_k\le h_k\le A_2f_k. $$
Then for NOR, Theorem  \ref{thm2.1} (i), (ii), (v) hold.

Indeed,  in  the proof of Theorem  \ref{thm2.1} (v), we only used
Condition (3) of Property (P). Hence for NOR, Theorem \ref{thm2.1}
(v) holds. In  the proofs of Theorem \ref{thm2.1} (i) and (ii), we
used Condition (3) of Property (P) and the monotonicity of
$\{h_k\}$. If Condition (2) is replaced by Condition $(2)'$, then
the inequalities \eqref{3.2-1} and \eqref{3.7} can be replaced by
the following inequalities $$\ln2 A_1^xd f_d^x\le \ln2\,
A_1^x\sum_{k=1}^d f_k^x \leq \ln2\sum\limits_{k=1}^d h_k^x\\
\le\ln\Big(\sum\limits_{j=1}^\infty \lambda_{d,j}^x\Big),$$ and
$$\ln\Big(\sum\limits_{j=1}^\infty \lambda_{d,j}^x\Big) \le M_x\sum_{k=1}^\infty h_k^x\le
M_xA_2^x\sum_{k=1}^\infty f_k^x.
$$
Using the above inequalities and the methods in the proofs of
Theorem  \ref{thm2.1} (i) and (ii), noting that
$$\liminf\limits_ {d\to\infty}\frac{\ln h_d^{-1}}{\ln
d}=\liminf\limits_ {d\to\infty}\frac{\ln f_d^{-1}}{\ln d},$$we can
prove that for NOR, Theorem  \ref{thm2.1} (i) and (ii) hold.
\end{rem}

\section{Applications of Theorem \ref{thm2.1}}

Consider the approximation problem ${\rm APP}=\{{\rm
APP}_d\}_{d\in\Bbb N}$,
$${\rm APP}_d\,\,: \,\, H_{K_d}\to L_2([0,1]^d)\quad \text{with} \quad {\rm APP}_df=f,$$
where $H_{K_d}$ is a Hilbert space related to the kernels $K_d$
which  are of tensor product and  correspond to Euler kernels,
Wiener kernels, Korobov kernels, Gaussian kernels, and analytic
Korobov kernels. This section is devoted to giving the
applications of Theorem 2.1 to these cases.

\subsection{Function approximation with Euler
kernels}

\

In this subsection we consider  multivariate approximation
problems  with Euler  kernels. Assume that  ${\bf r}=\{r_k\}_{k\in
\Bbb N}$  is a sequence of nondecreasing nonnegative integers
satisfying
\begin{equation}\label{4.0}0\le r_1\le r_2\le r_3 \le \dots.\end{equation} Let $H(K^E_{d,{\bf r}})$ be the
reproducing kernel Hilbert space with reproducing kernel being the
Euler kernel
$$K^E_{d,{\bf r}}(\x,\y)=\prod\limits_{k=1}^d K^E_{1,r_k}(x_k,y_k), \ \ \x,\y \in[0,
1]^d,$$where $$K^E_{1,r}(x, y)=\int_{[0,1]^r}\min(x, s_1)\min(s_1,
s_2)\dots\min(s_r, y)ds_1ds_2\dots ds_r, \ \ r\in \Bbb N_0 $$ is
the Euler kernel (see \cite[pp. 222-226]{NW3}). If $r=0$, we  get
the standard Wiener kernel
$$K^E_{1,0}(x, y)= \min(x,y).$$
If $X_r^E(t)$ is  a Gaussian random process with zero mean and
covariance kernel $K_{1,r}^E$, then it is called
  the univariate integrated Euler
process.

The reproducing kernel Hilbert space $H(K^E_{1,{ r}})$ is equal to
the space of functions $f: [0,1]\to \Bbb R$ such that the $r$th
derivative of $f$ is absolutely continuous and the $(r+1)$st
derivative of $f$ belongs to $L_2[0,1]$ and $f$ satisfies the
following boundary conditions
$$f(0)=f'(1)=f''(0)=\cdots=f^{(r)}(s_r)=0,$$
where $s_r=0$ if $r$ is even and $s_r=1$ if $r$ is odd. The inner
product of $H(K^E_{1,{ r}})$ is given by
$$\langle f,g\rangle_{r}=\int_0^1 f^{(r+1)}(x)g^{(r+1)}(x)dx \ \ {\rm for \ \ all} \ \ f,g\in H(K^E_{1,{ r}}).$$

Obviously, the space $H(K^E_{d,{\bf r}})$ is a tensor product of
the univariate $H(K^E_{1,{ r_j}})$, i.e.,
$$H(K^E_{d,{\bf r}})=H(K^E_{1,r_1})\otimes H(K^E_{1,r_2})\otimes\cdots\otimes H(K^E_{1,r_d}).$$

 We consider the multivariate approximation problem
${\rm APP}=\{{\rm APP}_d\}_{d\in \Bbb N}$ which is defined via the
embedding operator
$$ {\rm APP}_d: H(K^E_{d,{\bf r}})\to L_{2}([0,1]^d)\ \ {\rm with}\ \  {\rm APP}_d\,
f=f.$$

For the above approximation problem ${\rm APP}=\{{\rm APP}_d\}$,
the eigenvalues of the  operator $W_d={\rm APP}_d^*\ {\rm APP}_d$
  are given by (see \cite{GHT} or \cite[pp. 222-226]{NW3})
$$\big\{\lz^E_{d,\j} \big\}_{\j\in \Bbb N^d}=\big\{\lz^E(1, j_1)\lz^E(2, j_2)\dots\lz^E(d, j_d) \big\}_{(j_1,\dots, j_d )\in\Bbb N^d},$$
where
\begin{equation}\label{4.1}\lz^E(k,j)=\bigg[\frac{1}{\pi(j-\frac{1}{2})}\bigg]^{2r_k+2},\quad j\in \Bbb N
.\end{equation}

Now we verify that the above approximation problem APP satisfies
Conditions (2) and (3) of Property (P). First we note that the
sequence  $\{h^E_k\},\ h^E_k=\frac {\lz^E(k,2)}{
\lz^E(k,1)}=\frac{1}{3^{2r_k+2}}$ is nonincreasing due to
\eqref{4.0}. Next we have for $x>\frac{1}{2r_1+2}$,
\begin{equation*}\begin{aligned}&\sup_{k\in  \Bbb N}H^E(k,x)=\sup_{k\in  \Bbb N}\sum\limits_{j=2}^\infty \Big(\frac{\lambda^E(k,j)}{\lambda^E(k,2)}\Big)^x
\\&=\sup_{k\in  \Bbb N}\sum\limits_{j=2}^\infty\Big(\frac{3}{2j-1}\Big)^{x(2r_k+2)}
\leq\sum\limits_{j=2}^\infty\Big(\frac{3}{2j-1}\Big)^{2x(r_1+1)}<\infty.\end{aligned}\end{equation*}
It follows that APP satisfies Condition (3) of Property (P) with
$h^E_k=\frac{1}{3^{2r_k+2}}$ and  $\tau_0=\frac{1}{2r_1+2}$. By
Remark 2.2, we have the following corollary.

\begin{cor}\label{thm4.0} Consider the   approximation problem $ {\rm APP}=\big\{{\rm
APP}_d\big\}_{d\in\Bbb N}$ in the space $H(K^E_{d,{\bf r}})$ with
sequence ${\bf r}=\{r_k\}_{k\in \Bbb N}$ satisfying \eqref{4.0}.
Then for the normalized error criterion, we have

(1) ${\rm APP}$ is strongly polynomially tractable iff
\begin{equation*}A_*=\liminf\limits_{k\to\infty}\frac{(2r_k+2)\ln3}{\ln k}=2\ln 3\,\liminf\limits_{k\to\infty}\frac{r_k}{\ln k}>0,\end{equation*}
 and the exponent of SPT is
 $$ p^*=\max\big\{\frac{2}{A_*},\frac 1{r_1+1}\big\}.$$

(2) ${\rm APP}$ is strongly polynomially tractable iff it is
polynomially tractable.

\vskip 2mm

(3) ${\rm APP}$ is quasi-polynomially tractable for   all
sequences ${\bf r}=\{r_k\}_{k\in \Bbb N}$ satisfying \eqref{4.0}.
Furthermore, the exponent of QPT is
$$ t^*=\max\big\{\frac{2}{\lim\limits_{k\to\infty}(2r_k+2)\ln3},\frac 1{r_1+1}\big\}=\frac 1{r_1+1}.$$

(4) Obviously, QPT implies UWT,  WT,  and $(s,t)$-WT for any
positive $s$ and $t$, for all sequences ${\bf r}=\{r_k\}_{k\in
\Bbb N}$ satisfying \eqref{4.0}.

\end{cor}

\begin{rem}It is easy to see that for $x\ge\frac1{r_1+1}$,
$$\sup_{d\in\Bbb N}\sum_{j=1}^\infty (\lz_{d,j}^E)^x=\sup_{d\in\Bbb N}\prod_{k=1}^d
 \sum_{j=1}^\infty \bigg[\frac{1}{\pi(j-\frac{1}{2})}\bigg]^{(2r_k+2)x}\le \sup_{d\in\Bbb N}\prod_{k=1}^d
 \sum_{j=1}^\infty \bigg[\frac{1}{\pi(j-\frac{1}{2})}\bigg]^{2} <1. $$
Here, we use the equality that $$\sum_{j=1}^\infty
\frac1{(2j-1)^2}=\frac{\pi^2}8.$$ It follows from \cite[Theorem
5.1]{NW1} that for ABS, APP is strongly polynomially tractable for
all sequences ${\bf r}=\{r_k\}_{k\in \Bbb N}$ satisfying
\eqref{4.0}.

 Next we consider the exponent of SPT. Set $$G(x)=\sum_{j=1}^\infty
 \Big[\frac{1}{\pi(j-\frac{1}{2})}\Big]^{x},\ \  \overline
 r=\lim\limits_{k\to\infty}r_k.$$Then $G(x)$ is
 decreasing on $(1,+\infty)$, and $G(2)=1/2$. So there exists a
 unique
 $\xi_0\in(1,2)$ such that $G(\xi_0)=1$.  Then
 $$\sup_{d\in\Bbb N}\sum_{j=1}^\infty (\lz_{d,j}^E)^x=\sup_{d\in\Bbb N}\prod_{k=1}^d
G((2r_k+2)x) \ \ \left\{\begin{matrix}
 & <\infty, \qquad\text{if } x>\max\big\{\frac{\xi_0}{2\overline r+2}, \frac1{2r_1+2}\big\},\\
 &=\infty, \qquad\text{if } x<\max\big\{\frac{\xi_0}{2\overline r+2}, \frac1{2r_1+2}\big\}.
\end{matrix}\right. $$ This means that the exponent of SPT is
$$\max\big\{\frac{\xi_0}{\overline r+1},\frac1{r_1+1}\big\}.$$
\end{rem}

\subsection{Function approximation with Wiener
kernels}

\

In this subsection we consider  multivariate approximation
problems  with Wiener kernels.
 Let $H(K^W_{d,{\bf r}})$ be the
reproducing kernel Hilbert space with reproducing kernel being the
Wiener kernel
$$K^W_{d,{\bf r}}(\x,\y)=\prod\limits_{k=1}^d K^W_{1,r_k}(x_k,y_k), \ \ \x,\y \in[0,
1]^d,$$ where $$K^W_{1,r}(x,
y)=\int_0^{\min(x,y)}\frac{(x-u)^r}{r!}\frac{(y-u)^r}{r!}
du=\int_0^1\frac{(x-u)^r_+}{r!}\frac{(y-u)^r_+}{r!} du$$ is the
Wiener kernel, $x,y\in [0,1]$, $t_+=\max\{t,0\}$, and ${\bf
r}=\{r_k\}_{k\in \Bbb N}$  is a sequence of nondecreasing
nonnegative integers satisfying \eqref{4.0}.  If $r=0$, we also
get the standard Wiener kernel
$$K^W_{1,0}(x, y)= \min(x,y).$$If $X_r^W(t)$ is  a Gaussian random process with zero mean and
covariance kernel $K_{1,r}^W$, then it is called
  the univariate integrated Wiener
process. Integrated Euler and Wiener processes have many
applications in probability theory, statistics, physics.

The reproducing kernel Hilbert space $H(K^W_{1,{ r}})$ is equal to
the space of functions $f: [0,1]\to \Bbb R$ such
that $f$ is absolutely continuous and satisfies the following boundary conditions
$$f(0)=f'(0)=f''(0)=\cdots=f^{(r)}(0)=0.$$
The inner product of $H(K^W_{1,{ r}})$ is given by
$$\langle f,g\rangle_{r}=\int_0^1 f^{(r+1)}(x)g^{(r+1)}(x)dx \ \ {\rm for \ \ all} \ \ f,g\in H(K^W_{1,{ r}}).$$

Obviously, the space $H(K^W_{d,{\bf r}})$ is a tensor product of
the univariate $H(K^W_{1,r_j})$, i.e.,
$$H(K^W_{d,{\bf r}})=H(K^W_{1,r_1})\otimes H(K^W_{1,r_2})\otimes\cdots\otimes H(K^W_{1,r_d}).$$

 We consider the multivariate approximation problem
${\rm APP}=\{{\rm APP}_d\}_{d\in \Bbb N}$ which is defined via the
embedding operator
$$ {\rm APP}_d: H(K^W_{d,{\bf r}})\to L_{2}([0,1]^d)\ \ {\rm with}\ \  {\rm APP}_d\,
f=f.$$

For the above approximation problem ${\rm APP}=\{{\rm APP}_d\}$,
the eigenvalues of the  operator $W_d={\rm APP}_d^*\ {\rm APP}_d$
  are given by (see \cite{LPW1, LPW2})
$$\big\{\lz^W_{d,\j} \big\}_{\j\in \Bbb N^d}=\big\{\lz^W(1, j_1)\lz^W(2, j_2)\dots\lz^W(d, j_d) \big\}_{(j_1,\dots, j_d )\in\Bbb N^d},$$
where
\begin{equation*}\lz^W(k,j)=\bigg[\frac{1}{\pi(j-\frac{1}{2})}\bigg]^{2r_k+2}+\mathcal{O}(j^{-(2r_k+3)})
, \ \ j\to \infty,
\end{equation*} and for two positive sequences $f,g:\Bbb N\to [0,\infty),$
$$f(k)=\mathcal{O}(g(k)),\ \ k\to \infty$$
means that there exists two constants $C>0$ and $k_0\in \Bbb N$
for which $f(k)\leq Cg(k)$ holds for any $k\geq k_0$,
$$f(k)=\Theta(g(k)),\ \ k\to \infty$$ mean that
$$ f(k)=\mathcal{O}(g(k))\ \ {\rm and} \ \ g(k)=\mathcal{O}(f(k)),\ \ k\to \infty$$

Now we verify that the above approximation problem APP satisfies
Conditions (2) and (3) of Property (P). It was proved in
\cite{LPW2} that the  sequence $ \{h^W_k\}$ satisfies $$
h^W_k=\frac {\lz^W(k,2)}{
\lz^W(k,1)}=\Theta(r_k^{-2})=\Theta((1+r_k)^{-2}),\ \ k\to
\infty.$$ We conclude that the approximation problem APP satisfies
Condition $(2)'$ defined in Remark 3.1 with $f_k^W=(1+r_k)^{-2}, \
\ k\in\Bbb N.$

It follows from  \cite[Thm. 4.1]{LPW2} that for $x\in(3/5,1]$,
$$A_x:=\sup_{k\in\Bbb
N}\sum_{j=3}^\infty\Big(\frac{\lz^W(k,j)}{\lz^W(k,2)}\Big)^x<\infty.$$
We conclude that for $x>\frac{3}{5}$
\begin{equation*}\begin{aligned}&\sup_{k\in  \Bbb N}H^W(k,x)=\sup_{k\in  \Bbb N}\sum\limits_{j=2}^\infty \Big(\frac{\lambda^W(k,j)}{\lambda^W(k,2)}\Big)^x
\\&=1+\sup_{k\in  \Bbb N}\sum\limits_{j=3}^\infty\Big(\frac{\lambda^W(k,j)}{\lambda^W(k,2)}\Big)^x=1+A_x
<\infty.\end{aligned}\end{equation*} It follows that APP satisfies
Condition (3) of Property (P) with $\tau_0\in [0,\frac{3}{5}]$. By
Remark 3.1, we have the following corollary.

\begin{cor}\label{thm4.6} Consider the   approximation problem $ {\rm APP}=\big\{{\rm
APP}_d\big\}_{d\in\Bbb N}$ in the space $H(K^W_{d,{\bf r}})$ with
sequence ${\bf r}=\{r_k\}_{k\in \Bbb N}$ satisfying \eqref{4.0}.
Then for the normalized error criterion, we have

(1) ${\rm APP}$ is strongly polynomially tractable iff
\begin{equation*}A_*=\liminf\limits_{k\to\infty}\frac{2\ln(1+r_k)}{\ln k}=2\liminf\limits_{k\to\infty}\frac{\ln r_k}{\ln k}>0,\end{equation*}
and the exponent of SPT is  $$ p^*=\max\big\{\frac{2}{A_*},2{\tau_0}\big\}.$$

(2) ${\rm APP}$ is strongly polynomially tractable iff it is
polynomially tractable.

\vskip 2mm

(3) ${\rm APP}$ is $(s,t)$-weakly tractable with $s>0$ and $t>1$
for all sequences ${\bf r}=\{r_k\}_{k\in \Bbb N}$ satisfying
\eqref{4.0}.
\end{cor}

\begin{rem}We do not know the exact value of $\tau_0$. It is open. \end{rem}

\subsection{Function approximation with Korobov
kernels}

\

In this subsection we consider  multivariate approximation
problems  with Korobov  kernels. First we recall  definition of
Korobov spaces (see \cite[Appendix A]{NW1}). Let ${\bf
r}=\{r_k\}_{k\in \Bbb N} $ and ${\bf g}=\{g_k\}_{k\in \Bbb N}$ be
two sequences  satisfying
\begin{equation}\label{4.11} 1 \geq g_1\geq g_2\geq \cdots\geq g_k \geq \cdots > 0,\end{equation}
and
\begin{equation}\label{4.12} 0<r_1\le r_2\le \cdots\le r_k\le \cdots. \end{equation}
 For $d=1,2,\cdots$, we define the spaces
 $$ H_{d,{\bf r,g}}=H_{1,r_1,g_1}\otimes H_{1,r_2,g_2}\otimes\cdots\otimes H_{1,r_d,g_d}.$$
 Here $H_{1,\az,\beta}$ is the Korobov space of univariate complex valued functions $f$ defined on $[0,1]$ such that
$$\|f\|_{H_{1,\alpha,\beta}}^2:=|\hat f(0)|^2+\beta^{-1}\sum\limits_{h\in \Bbb Z,h\neq0}|h|^{2\az}|\hat f(h)|^2<\infty , $$where $\beta\in(0, 1]$ is a scaling parameter,
and $\alpha>0$ is a smoothness parameter,
$$\hat f(h)=\int_0^1\exp (-2\pi ihx)f(x)dx\ \  {\rm for}\ \ h\in \Bbb
Z$$are the  Fourier coefficients of $f$, $i=\sqrt{-1}$.  If $\az>
\frac{1}{2}, $ then $H_{1,\az,\beta}$ consists of $1$-periodic
functions  and  is a reproducing kernel Hilbert space with
reproducing kernel
$$R_{\alpha,\beta}(x,y):=1+2\beta\sum\limits_{j=1}^\infty j^{-2\alpha}\cos(2\pi j(x-y)),\  \ x,y\in [0,1].$$
 If $\az$ is an
integer, then $H_{1,\az,\beta}$ consists of $1$-periodic functions
$f$ such that $f^{(\az-1)}$ is absolutely continuous,  $f^{(\az)}$
belongs to $L_2([0,1])$, and
$$\|f\|_{H_{1,\alpha,\beta}}^2=\big|\int_{[0,1]}f(x)dx\big|^2+ (2\pi)^{2\az}\beta^{-1} \int_{[0,1]} |f^{(\az)}(x)|^2dx.$$

For $d\ge2$ and two sequences ${\bf r}=\{r_k\}$ and ${\bf
g}=\{g_k\}$ satisfying \eqref{4.11} and \eqref{4.12}, the space
$H_{d, {\bf \az,\beta}}$ is a Hilbert space with the inner product
$$ \langle f,g\rangle_{H_{d,{\bf r,g}}}=\sum\limits_{{\bf h}\in \Bbb Z^d} \rho_{d,{\bf r,g}}({\bf h})\hat f({\bf h}) \overline{\hat g({\bf h})},$$
where $$\rho_{d,{\bf r,g}}({\bf
h})=\prod\limits_{j=1}^d(\delta_{0,h_j}+g_j^{-1}(1-\delta_{0,h_j}))|h_j|^{2r_j},$$
$\delta_{i,j}=\Big\{\begin{array}{ll}1, \ &i=j,\\ 0, &i\neq
j,\end{array}$\ \ and
$$\hat f(\h)=\int_{[0,1]^d}\exp (-2\pi i\h \cdot \x)f(\x)d\x\ \  {\rm for}\ \ \h\in \Bbb
Z^d$$are the  Fourier coefficients of $f$, $\x \cdot
\y=x_1y_1+\dots+x_dy_d$. If $r_1> \frac{1}{2}, $ then $H_{d,{\bf
r,g}}$ consists of $1$-periodic functions on $[0,1]^d$ and  is a
reproducing kernel Hilbert space with  reproducing kernel
\begin{align*}K_{d, {\bf r,g}}(\x,\y)&=\prod_{k=1}^d
R_{r_k,g_k}(x_k,y_k)\\
&=\prod_{k=1}^d\Big(1+2g_k\sum\limits_{j=1}^\infty
j^{-2r_k}\cos(2\pi j(x_k-y_k))\Big),\
 \ \x,\y\in [0,1]^d.\end{align*}
 For integers $r_j$, the inner product of $H_{d, {\bf r,g}}$ can be expressed in terms of
 derivatives, see \cite[Appendix A]{NW1}. If $\{r_k\}$ is
 nondecreasing and $g_k=(2\pi)^{-2r_k}$, then we obtain the
 specific Korobov spaces $H_{d,{\bf r}}$ which were given in
 \cite{PW, S2}.

 In this subsection, we consider the multivariate approximation problem
${\rm APP}=\{{\rm APP}_d\}_{d\in \Bbb N}$ which is defined via the
embedding operator
$$ {\rm APP}_d: H_{d,{\bf r,g}}\to L_{2}([0,1]^d)\ \ {\rm with}\ \  {\rm APP}_d\,
f=f.$$

For the above approximation problem ${\rm APP}=\{{\rm APP}_d\}$,
the eigenvalues of the  operator $W_d={\rm APP}_d^*\ {\rm APP}_d$
  are given by (see \cite[p. 184]{NW1})
$$\big\{\lz_{d,\j} \big\}_{\j\in \Bbb N^d}=\big\{\lz(1, j_1)\lz(2, j_2)\dots\lz(d, j_d) \big\}_{(j_1,\dots, j_d )\in\Bbb N^d},$$
where $\lz(k,1)=1$, and
$$\lz(k,2j)=\lz(k,2j+1)={g_k}\,{j^{-2r_k}}\in(0,1],\ \ j\in  \Bbb N,.$$

Now we verify that the above approximation problem APP satisfies
Property (P). First we note that $\lz(k,1)=1$, and the sequence
$\{h_k\},\ h_k=\frac {\lz(k,2)}{ \lz(k,1)}=g_k$ is nonincreasing
due to \eqref{4.11}. Next we have for $x>\frac 1{2r_1}$,
\begin{equation*}\begin{aligned}\sup_{k\in  \Bbb N}H(k,x)&
=\sup_{k\in  \Bbb N}\sum\limits_{j=2}^\infty
\Big(\frac{\lambda(k,j)}{\lambda(k,2)}\Big)^x
 =2\sup_{k\in  \Bbb N}\sum\limits_{j=1}^\infty j^{-2r_kx}
 \\&=H(1,x)= 2\sum\limits_{j=1}^\infty j^{-2r_1 x}
 =2\zeta(2r_1x)<\infty,\end{aligned}\end{equation*}
where $\zeta(s)=\sum_{j=1}^\infty j^{-s}$ is the Riemann zeta
function which is well-defined only for $s> 1$. This means that
the above approximation problem APP has Property (P) with
$h_k=g_k$ and $\tau_0=\frac 1{2r_1}$.
 By Theorem 2.1, we have the following corollary.

\begin{cor}\label{cor4.1}Consider the above   approximation problem $ {\rm APP}=\big\{{\rm
APP}_d\big\}_{d\in\Bbb N}$ in the worst case setting with the
sequences ${\bf r}$ and ${\bf g}$
 satisfying \eqref{4.11} and \eqref{4.12}. Then
for the absolute error criterion or the normalized error
criterion, we have

(1) ${\rm APP}$ is strongly polynomially tractable iff
\begin{equation*} A_*=\liminf\limits_ {k\to\infty}\frac{\ln g_k^{-1}}{\ln k}>0,\end{equation*}
and the exponent of SPT is \begin{equation*} p^*=\max\big\{\frac{2}{A_*},\frac{1}{r_1}\big\}.\end{equation*}

(2) ${\rm APP}$ is strongly polynomially tractable iff it is
polynomially tractable.

\vskip 2mm

(3) ${\rm APP}$ is  quasi-polynomially tractable iff
$B:=\lim\limits_{k\to\infty}\ln g_k^{-1}>0$ iff $g_k\not\equiv 1$.
Furthermore, the exponent of QPT is
$$ t^*=\max\big\{\frac{2}{B},\frac{1}{r_1}\big\}.$$

(3) ${\rm APP}$ is  quasi-polynomially tractable iff ${\rm APP}$
is  uniformly weakly tractable, iff ${\rm APP}$ is $(s,t)$-weakly
tractable with $s>0$ and $t\in([0,1]$, and iff  ${\rm APP}$ is
weakly tractable.

\vskip2mm

(4) ${\rm APP}$ is $(s,t)$-weakly tractable with $s>0$ and $t>1$
for all sequences ${\bf r}$ and ${\bf g}$
 satisfying \eqref{4.11} and \eqref{4.12}.

 \vskip2mm

(5) If $g_k\equiv1$, then ${\rm APP}$ suffers from the curse of
dimensionality.
\end{cor}

\begin{rem}In \cite{PW, S2}, the authors investigated tractability of the approximation problem
$$ {\rm APP}_d: H_{d,{\bf r}}\to L_{2}([0,1]^d)\ \ {\rm with}\ \  {\rm APP}_d\,
f=f,$$ where $H_{d,{\bf r}}=H_{d,{\bf r,g}}$, ${\bf r}=\{r_k\}$
and ${\bf g}=\{g_k\}$ satisfy $0<r_1\le r_2\le \dots,$
$g_k=(2\pi)^{-2r_k}, k=1,2,\dots$. They obtained the following
results.

For the absolute error criterion or the normalized  error
criterion, we have

\vskip 2mm

 $\bullet$ SPT holds iff PT  holds iff
$$R=\limsup\limits_{k\to\infty}\frac{\ln k}{r_k}<+\infty,$$ with the
exponent of SPT $$p^*=\max \big\{\frac 1{r_1}, \frac{R}{\ln 2\pi}\big\}.$$

\vskip 2mm

$\bullet$ QPT holds iff UWT holds iff WT holds iff $r_1>0$.
Furthermore, the exponent of QPT is
$$t^*=\frac 1{r_1}.$$
\vskip 2mm

We remark that the above results follows directly from Corollary
\ref{cor4.1}. Hence Corollary \ref{cor4.1} is a generalization of
the above conclusion.
  \end{rem}

\subsection{Function approximation  with Gaussian kernels}

\

In this subsection we consider  multivariate approximation
problems  with Gaussian kernels. Let $H(K_{d,\gaa})$ be the
reproducing kernel Hilbert space with the Gaussian kernel
\begin{equation*}
K_{d,\gaa}(\x,\y)=\prod_{j=1}^{d}K_{1,\gz_j}(x_j, y_j), \ \x,\y\in
\R^d,
\end{equation*}where $$K_{1,\gz}(x, y)=\exp (-\ga^2(x-y)^2), \ \ \gz>0,\ x,y\in
\R,$$ and $\gaa=\{\ga_j^2\}_{j\in \Bbb N}$ is a given sequence of
shape parameters  not depending on $d$ and satisfying
\begin{equation}\label{4.21}
\ga_1^2\geq \ga_2^2\geq\dots>0.
\end{equation}Obviously, the space $H(K_{d,\gaa})$ is a tensor product of the univariate $H(K_{1,\gamma_j})$, i.e.,
$$H(K_{d,\gaa})=H(K_{1,\gamma_1})\otimes H(K_{1,\gamma_2})\otimes\cdots\otimes H(K_{1,\gamma_d}).$$ The reproducing kernel Hilbert space
$H(K_{d,\gaa})$ has been widely used in many fields such as
numerical computation, statistical learning, and engineering (see
e.g., \cite{FHW1, FSK, RW}).

 We consider the multivariate approximation problem
${\rm APP}=\{{\rm APP}_d\}_{d\in \Bbb N}$ which is defined via the
embedding operator
$$ {\rm APP}_d: H(K_{d,\gaa})\to L_{2,d}\ \ {\rm with}\ \  {\rm APP}_d\,
f=f,$$where  $$L_{2,d}=\Big\{ f\ \big| \
\|f\|_{L_{2,d}}=\bigg(\int_{\R^d} |f(\x)|^2\prod_{j=1}^d \frac
{\exp {(-x_j^2)}}{\sqrt{\pi}}\,d\x\bigg)^{1/2}<\infty\Big\}$$ is a
separable Hilbert space of real-valued functions on $\R^d$ with
inner product $$\langle f,g\rangle_{L_{2,d}}= \int_{\R^d} f(\x)
g(\x)\prod_{j=1}^d \frac {\exp {(-x_j^2)}}{\sqrt{\pi}}\,d\x.$$

For the above approximation problem ${\rm APP}=\{{\rm APP}_d\}$,
the eigenvalues of the  operator $W_d={\rm APP}_d^*\ {\rm APP}_d$
  are given by (see \cite{FHW1, RW, SW})
  $$\big\{\lz_{d,\j} \big\}_{\j\in \Bbb N^d}=\big\{\lz(1, j_1)\lz(2, j_2)\dots\lz(d, j_d) \big\}_{(j_1,\dots, j_d )\in\Bbb N^d},$$
where  \begin{equation}\label{4.22}\lz
(k,j_k)=(1-\oz_{\ga_k})\oz_{\ga_k}^{j_k-1},\quad {\rm and} \quad
\oz_\ga=\frac{2\ga^2}{1+2\ga^2+\sqrt{1+4\ga^2}}\ \ {\rm for}\
\ga>0.
\end{equation} Clearly, $0<\oz_\ga<1$ for $\ga>0$ and  $\oz_\ga$ is
an increasing function of $\ga$ and $\oz_\ga$ tends to $0$ iff
$\ga$ tends to $0$.  We also have
\begin{equation}\label{4.23}\lim\limits_{\gz\to0} \frac {\oz_\ga}{\gz^2}=1, \ \ {\rm and}\ \ \lim_{\gz\to0} \frac
{\ln \oz_\ga^{-1}}{\ln \gz^{-2}}=1.\end{equation}

Now we verify that the above approximation problem APP satisfies
Conditions (2) and (3) of Property (P). First we note that the
sequence  $\{h_k\},\ h_k=\frac {\lz(k,2)}{ \lz(k,1)}=\oz_{\ga_k}$
is nonincreasing due to \eqref{4.21} and \eqref{4.22}. Next we
have for any $x>0$,
$$H(k,x)=\sum\limits_{j=2}^\infty \Big(\frac{\lambda(k,j)}{\lambda(k,2)}\Big)^x
=\sum\limits_{j=2}^\infty
\Big(\frac{\oz_{\ga_k}^{j-1}}{\oz_{\ga_k}}\Big)^x
=\sum\limits_{j=2}^\infty
\oz_{\ga_k}^{x(j-2)}=\frac1{1-\oz_{\ga_k}^{x}},$$ so that for
$x>0$,
$$\sup_{k\in  \Bbb N}H(k,x)=H(1,x)= \frac1{1-\oz_{\ga_1}^{x}}<\infty.$$It
follows that APP satisfies Condition (3) of Property (P) with
$h_k=\oz_{\ga_k}$ and  $\tau_0=0$. By Remark 2.2 and \eqref{4.23},
we have the following corollary.

\begin{cor}\label{thm4.4} Consider the   approximation problem $ {\rm APP}=\big\{{\rm
APP}_d\big\}_{d\in\Bbb N}$ in the space $H(K_{d,\gaa})$ with shape
parameters $\gaa=\{\gamma_j\}$ satisfying \eqref{4.21}. Then for
the normalized error criterion, we have

(1) ${\rm APP}$ is strongly polynomially tractable iff
\begin{equation*}r(\gaa)=\liminf\limits_{k\to\infty}\frac{\ln\oz_{\ga_k}^{-1}}{\ln k}= \liminf\limits_{k\to\infty}\frac{\ln \ga_k^{-2}}{\ln k}>0,\end{equation*}
 and the exponent of SPT is $\frac2{r(\gaa)}$.

\vskip 2mm

 (2) ${\rm APP}$ is strongly polynomially tractable iff
it is polynomially tractable.

\vskip 2mm

(3) ${\rm APP}$ is quasi-polynomially tractable for   all shape
parameters.  Furthermore, the exponent of QPT is
$$ t^*=\frac{2}{\lim\limits_{k\to\infty}\ln\frac{1+2\ga_k^2+\sqrt{1+4\ga_k^2}}{2\ga_k^2}}.$$

(4) Obviously, ${\rm APP}$ is  quasi-polynomially tractable
implies $(s,t)$-weakly tractable for any positive $s$ and $t$, as
well as uniformly weakly tractable  and   weakly tractable for all
shape parameters.\end{cor}

\begin{rem} In Corollary \ref{thm4.4} we obtain the exact value of the exponent of QPT. This result is new. The other results in Corollary \ref{thm4.4} were obtained
in \cite{FHW1}. Indeed, the authors in  \cite{FHW1} obtained
complete results about the tractability of APP as follows.

\vskip 3mm

For the absolute error criterion, we have

\vskip 2mm

 $\bullet$  SPT holds for all shape parameters with the
exponent $\min\big\{2, \frac{2}{r(\gaa)}\big\},$ where $r(\gaa)$ is
defined by
\begin{align*}r(\gaa)&=\sup\,\big\{\dz>0\ | \ \sum_{j=1}^\infty
\ga_{j}^{2/\dz}<\infty\big\}\\ &=\sup\,\big\{\beta\ge 0\ |\
\lim_{j\to\infty}j^\beta\gz_j^2=0\big\}=\liminf\limits_{j\to\infty}\frac{\ln
\ga_j^{-2}}{\ln j},
\end{align*}
where the last equality follows from \cite{CW1}. \vskip 2mm

$\bullet$ Obviously, SPT implies all PT, QPT, WT, $(s,t)$-WT for
any positive $s$ and $t$, as well as UWT and WT, for all shape
parameters.

\vskip 3mm

For the normalized error criterion, we have

\vskip 2mm

 $\bullet$ SPT holds iff PT  holds iff
$r(\gaa)>0$, with the exponent $\frac{2}{r(\gaa)}$.

\vskip 2mm

$\bullet$ QPT holds for all shape parameters with the exponent
$$t^*\le
\frac{2}{\ln\frac{1+2\ga_1^2+\sqrt{1+4\ga_1^2}}{2\ga_1^2}}.$$

\vskip 2mm

$\bullet$ Obviously, QPT implies $(s,t)$-WT for any positive $s$
and $t$, as well as UWT and WT for all shape parameters.
\end{rem}

\subsection{Function approximation  with
 analytic Korobov
kernels}

\

In this subsection we consider  multivariate approximation
problems  with analytic Korobov kernels. Let $H(K_{d,{\bf a,b}})$
be the reproducing kernel Hilbert space with the analytic Korobov
kernel
$$K_{d,{\bf a,b}}(\x,\y)=\prod\limits_{k=1}^dK_{1,a_k,b_k}(x_k,y_k), \ \ \x,\y \in[0,
1]^d,$$where ${\bf a}=\{a_k\}_{k\in \Bbb N}$ and ${\bf
b}=\{b_k\}_{k\in \Bbb N}$ are two sequences of positive weights
satisfying
\begin{equation}\label{4.30}0<a_1\leq a_2\leq \cdots \leq a_k \leq\cdots ,\ \ {\rm and}\ \  b_*:=\inf\limits_{k\in \Bbb N}b_k>0,\end{equation} $K_{1,a_k,b_k}$ are univariate analytic Korobov
kernels,
$$K_{1,a,b}(x,y)=\sum\limits_{h\in \Bbb Z}\omega^{a|h|^b}\exp(2\pi ih(x-y)),\ \ x, y\in[0,
1].$$ Here $\omega \in(0, 1)$ is a fixed number,
$i=\sqrt{-1},a,b>0$. Hence, we have
$$K_{d,{\bf a,b}}(\x,\y)=\sum\limits_{\h\in \Bbb Z}\omega_{\bf h} \exp(2\pi i\h(\x-\y)),\ \ \x, \y\in[0, 1]^d,$$
with
$$\omega_\h=\omega^{\sum\limits_{k=1}^da_k|h_k|^{b_k}},\ \forall\ \ \h=(h_1, h_2,\cdots, h_d)\in \Bbb Z^d,$$
for fixed $\omega \in (0, 1)$.

Obviously, the space $H(K_{d,{\bf a,b}})$ is a tensor product of
the univariate $H(K_{1,a_k,b_k})$, i.e.,
$$H(K_{d,{\bf a,b}})=H(K_{1,a_1,b_1})\otimes H(K_{1,a_2,b_2})\otimes\cdots\otimes H(K_{1,a_d,b_d}).$$ The reproducing kernel Hilbert space
 $H(K_{d,{\bf a,b}})$ has been widely used in the study of
tractability and exponential convergence-tractability (see
\cite{DKPW,  IKPW, KPW, LX1}).

 We consider the multivariate approximation problem
${\rm APP}=\{{\rm APP}_d\}_{d\in \Bbb N}$ which is defined via the
embedding operator
$$ {\rm APP}_d: H(K_{d,{\bf a,b}})\to L_{2}([0,1]^d)\ \ {\rm with}\ \  {\rm APP}_d\,
f=f.$$

For the above approximation problem ${\rm APP}=\{{\rm APP}_d\}$,
the eigenvalues of the  operator $W_d={\rm APP}_d^*\ {\rm APP}_d$
  are given by
$$\big\{\lz_{d,\j} \big\}_{\j\in \Bbb N^d}=\big\{\lz(1, j_1)\lz(2, j_2)\dots\lz(d, j_d) \big\}_{(j_1,\dots, j_d )\in\Bbb N^d},$$
where $\lz(k,1)=1$, and
$$\lz(k,2j)=\lz(k,2j+1)=\oz^{a_kj^{b_k}},\quad j\in\Bbb N.$$

Now we verify that the above approximation problem APP satisfies
Property (P). First we note that $\lz(k,1)=1$, and the sequence
$\{h_k\},\ h_k=\frac {\lz(k,2)}{ \lz(k,1)}=\oz^{a_k}$ is
nonincreasing due to \eqref{4.30}. Next we have for any $x>0$,
$$\sup_{k\in  \Bbb N}H(k,x)=\sup_{k\in  \Bbb N}\sum\limits_{j=2}^\infty \Big(\frac{\lambda(k,j)}{\lambda(k,2)}\Big)^x
=2\sup_{k\in  \Bbb N}\sum\limits_{j=1}^\infty
\omega^{xa_k(j^{b_k}-1)} \leq 2\sum\limits_{j=1}^\infty
\omega^{xa_1(j^{b_*}-1)}.$$ Since
$$\omega^{xa_1(j^{b_*}-1)}=j^{-\frac{xa_1(j^{b_*}-1)\ln\omega^{-1}}{\ln
j}} \ \ {\rm and} \ \  \lim\limits_{j\to
\infty}\frac{xa_1(j^{b_*}-1)\ln\omega^{-1}}{\ln j}=\infty,$$
we  conclude that for $x>0$,
\begin{equation*}\sup_{k\in  \Bbb N}H(k,x)\leq M_x:=2\sum\limits_{j=1}^\infty  \omega^{x a_1(j^{b_*}-1)}<\infty.\end{equation*}
This means that the above approximation problem APP has Property
(P) with  $h_k=\omega^{a_k}$ and $\tau_0=0$.
 By Theorem 2.1, we have the following corollary.

\begin{cor}\label{cor4.5}Consider the above   approximation problem $ {\rm APP}=\big\{{\rm
APP}_d\big\}_{d\in\Bbb N}$ in the space $H(K_{d,{\bf a,b}})$ with
the sequences ${\bf a}$ and ${\bf b}$
 satisfying \eqref{4.30}. Then
for the absolute error criterion or the normalized error
criterion, we have

(1) ${\rm APP}$ is strongly polynomially tractable iff
\begin{equation*}A_*=\liminf\limits_{k\to\infty}\frac{{a_k}\ln\omega^{-1}}{\ln k}=\ln\omega^{-1}\liminf\limits_{k\to\infty}\frac{{a_k}}{\ln k}>0,\end{equation*}
 and the exponent of SPT is
$$ p^*=\frac{2}{A_*}.$$

(2) ${\rm APP}$ is strongly polynomially tractable iff it is
polynomially tractable.

\vskip 2mm

(3) ${\rm APP}$ is quasi-polynomially tractable for all sequences
${\bf a}$ and ${\bf b}$ satisfying \eqref{4.30}.  Furthermore, the
exponent of QPT is
    $$t^*=\frac{2}{B},\ \  B=\ln \omega^{-1} \lim_{k\to\infty}a_k.$$

(4) Obviously, quasi-polynomially tractable implies $(s,t)$-weakly
tractable for any positive $s$ and $t$, as well as uniformly
weakly tractable and weakly tractable for  all sequences ${\bf a}$
and ${\bf b}$ satisfying \eqref{4.30}. \end{cor}

\begin{rem}Corollary \ref{cor4.5} is not new. All results in Corollary \ref{cor4.5} were obtained
in \cite[Theorem 3.1]{LX1}.  See \cite{IKPW, KPW, LX1} for
background and more information.
\end{rem}

\noindent{\bf Acknowledgment}  This work was supported by the
National Natural Science Foundation of China (Project no.
11671271) and
 the  Natural Science Foundation of Beijing Municipality (1172004).


\begin{thebibliography}{99}



\bibitem{CW1}J. Chen, H. Wang, Average case tractability of  multivariate  approximation with
Gaussian  kernels, J. Approx. Theory 239 (2019) 51-71.


\bibitem{CWZ}J. Chen, H. Wang, J. Zhang, Average case $(s, t)$-weak tractability of non-homogeneous tensor product
problems,  J. Complexity 49 (2018) 27-45.

\bibitem {DKPW}  J. Dick, P. Kritzer, F. Pillichshammer, H. Wo\'{z}niakowski, Approximation of analytic functions in
Korobov spaces, J. Complexity  30 (2014) 2-28.

\bibitem{FHW1}G. E. Fasshauer, F. J. Hickernell, H. Wo\'zniakowski, On dimension-independent rates of
convergence for function approximation with Gaussian kernels, SIAM
J. Numer. Anal. 50 (2012) 247-271.

\bibitem{FHW2} G. E. Fasshauer, F. J. Hickernell, H. Wo\'zniakowski, Average case approximation:
 convergence and tractability of Gaussian kernels,   Monte Carlo and Quasi-Monte Carlo 2010, eds. L. Plaskota and H. Wo\'zniakowski, Springer Verlag, 2012, 329-345.

\bibitem{FSK}A. I. J. Forrester, A. S\'obester, and A. J. Keane, Engineering Design via Surrogate Modelling: A Practical Guide, Wiley, Chichester, 2008.

\bibitem{GHT}F. Gao, J. Hanning, F. Torcaso, Integrated Brownian motions and exact $L_2$-small balls, Ann. Probab. 31 (2003) 1320-1337.
\bibitem{IKPW}C. Irrgeher, P. Kritzer, F. Pillichshammer, H. Wo\'zniakowski,
 Tractability of multivariate approximation defined over Hilbert spaces with exponential weights, J. Approx. Theory 207 (2016) 301-338.

\bibitem{KPW}P. Kritzer, F. Pillichshammer, H. Wo\'zniakowski, Tractability of multivariate analytic problems, in: P. Kritzer,
H. Niederreiter, F. Pillichshammer, A. Winterhof (Eds.), Uniform
Distribution and Quasi-Monte Carlo Methods. Discrepancy,
Integration and Applications, De Gruyter, Berlin, 2014, pp.
147-170.
\bibitem{LPW1}M. A. Lifshits, A. Papageorgiou, H. Wo\'zniakowski, Average case tractability of non-homogeneous tensor product problems, J. Complexity 28 (2012) 539-561.
\bibitem{LPW2}M. A. Lifshits, A. Papageorgiou, H. Wo\'zniakowski, Tractability of multi-parametric Euler and Wiener integrated processes, Probab. Math. Statist. 32 (2012) 131-165.

\bibitem{LX1} Y. Liu, G. Xu,  A note on tractability of multivariate analytic
problems, J. Comlexity 34 (2016) 42-49.
\bibitem{LX} Y. Liu, G. Xu, Average case tractability of a multivariate approximation problem, J. Comlexity, 43 (2017) 76-102.
\bibitem{NW1} E. Novak, H. Wo\'zniakowski, Tractablity  of Multivariate Problems, Volume I: Linear Information, EMS, Z\"urich, 2008.
\bibitem{NW2} E. Novak, H. Wo\'zniakowski, Tractablity  of Multivariate Problems, Volume II: Standard Information for Functionals, EMS, Z\"urich, 2010.
\bibitem{NW3} E. Novak, H. Wo\'zniakowski, Tractablity  of Multivariate Problems, Volume III: Standard Information for Operators, EMS, Z\"urich, 2012.

\bibitem{PW}A. Papageorgiou, H. Wo\'zniakowski, Tractability through increasing
smoothness, J. Complexity  26 (2010) 409-421.

\bibitem{RW}C. E. Rasmussen, C. Williams, Gaussian Processes for Machine Learning, MIT Press, 2006 (online version at http://www.gaussianprocess.org/gpml/).

\bibitem{S1} P. Siedlecki, Uniform weak tractability, J. Complexity 29 (6) (2013) 438-453.
\bibitem{S2}P. Siedlecki, Uniform weak tractability of multivariate problems with increasing smoothness, J. Complexity 30 (2014) 716-734.
\bibitem{S3}P. Siedlecki, $(s,t)$-weak tractability of Euler and Wiener integrated processes, J. Complexity  45 (2018) 55-66.
\bibitem{SiW}P. Siedlecki,  M. Weimar, Notes on $(s, t)$-weak tractability: a refined classification of problems with
    (sub)exponential information complexity, J. Approx. Theory 200 (2015) 227-258.
\bibitem{SW} I. H. Sloan, H. Wo\'zniakowski, Multivariate approximation for analytic functions with Gaussian kernels, J. Complexity 45
(2018) 1-21.

\bibitem {W}  H. Wo\'{z}niakowski, Tractablity and strong tractability of multivariate tensor product Problems,
J. Comput. Inform. 4 (1994) 1-19.

\bibitem{X1}G.  Xu, Quasi-polynomial tractability of linear problems in the average case setting, J. Complexity 30 (2014)
54-68.
\bibitem{X2}G.  Xu, Tractability of linear problems defined over Hilbert spaces, J. Complexity 30 (2014) 735-749.


\end{thebibliography}
\end{document}